\documentclass{article}

\usepackage{amsfonts}
\usepackage{amsmath}
\usepackage{graphicx}%

\numberwithin{equation}{section}

\newtheorem{theorem}{Theorem}[section]

\newtheorem{corollary}[theorem]{Corollary}
\newtheorem{lemma}[theorem]{Lemma}
\newtheorem{remark}[theorem]{Remark}

\newenvironment{proof}[1][Proof]{\noindent\textbf{#1.} }{\ \rule{0.5em}{0.5em}}

\newcommand{\ra}{\rightarrow}

\newcommand{\dist}{{\rm dist}}

\newcommand{\C}{{\rm \bf C}}
\newcommand{\R}{{\rm \bf R}}
\newcommand{\Rn}{{\R^n}}
\newcommand{\Rnp}{{\R^{n+1}_+}}

\begin{document}

\title{Classes of measures generated by capacities}
\author{G. Dafni\footnote{Supported  in part by NSERC and the Centre de
Recerca Matematica, Barcelona.},  G.E. Karadzhov\footnote{Supported
in part by Memorial and Concordia Universities, NSERC, and the
Centre de Recherches Mathematiques, Montreal.} \ and J.
Xiao\footnote{Supported in part by NSERC of Canada and Dean of
Science Start-up Funds of MUN (Canada).}}

\maketitle

\begin{abstract}
We introduce classes of measures in the half-space
$\mathbf{R}^{n+1}_+,$ generated by Riesz, or Bessel, or Besov
capacities in $\mathbf{R}^n$, and give a geometric
characterization as Carleson-type measures.
\end{abstract}

\section{Introduction}

Recall that a Carleson measure is a positive Borel measure $\mu$ on the
upper half-space $\Rnp = \Rn \times (0,\infty)$ satisfying
\begin{equation}
\label{Carleson}
\mu(TB) \leq C |B|
\end{equation}
for some constant $C < \infty$ and all balls $B$ in $\Rn$. Here $|B|$ is
the
$n$-dimensional Lebesgue measure of $B$ and $TB$ is the ``tent" over B:
$$TB = \{(x,t): x \in B, t <\dist(x,\partial B)\}.$$

Note that in condition (\ref{Carleson}), the balls may be replaced by
cubes.  From this, using a  Whitney decomposition and the additivity of
Lebesgue measure,  we can replace balls by any open sets (see \cite{S2},
Chapter II, Section 2.3.2).

A positive Borel measure $\mu$ on the upper half-space $\Rnp$ is called
a $\beta$-Carleson measure if
\begin{equation}
\label{Cbeta}
\mu(T(B)) \leq C |B|^\beta \quad
\end{equation}
for some constant $C < \infty$ and all balls $B$.
We denote the class of such measures by  $C^\beta$.   These classes
were considered in
\cite{DX}, \cite{KX}

When $\beta \geq 1$, we can again replace balls by open sets
since we have
$$\sum |B_k|^\beta \leq (\sum |B_k|)^\beta.$$

However, in the case $\beta < 1$, the condition on open sets
is stronger.  For example, in the upper-half plane, consider
the measure
$$\mu = \sum_{k = 1}^\infty \frac{1}{k} \delta_{(k,
\frac{1}{k^{1/\beta}})},$$
where $\delta_P$ denotes the unit mass at point $P$.
Then $\mu(T(I)) \leq C |I|^\beta$ for every interval
$I$
but one can construct open sets $O_K = \cup I_k$
with
$$|O_K| =  \sum_{k = 1}^K |I_k| \approx \sum_{k = 1}^K
\frac{1}{k^{1/\beta}} \leq C < \infty
$$ and
$$\mu(T(O_K)) \approx \sum_{k = 1}^K \frac{1}{k} \ra \infty.$$
The spaces $V^\beta$ of measures satisfying the condition (\ref{Cbeta})
for all bounded open sets were considered in \cite{AB},
\cite{AM}, \cite{BJ},\cite{J}.

We are interested in whether we can modify (\ref{Cbeta}) so as to
give a condition on open sets which is equivalent to the $\beta$-Carleson
condition on balls.  An example of such a result is the following:

\begin{lemma}[DX]
\label{lemma-Cbeta}
For $0 < \beta \leq 1$,  $\mu \in C^\beta$ if and
only if there exists $C < \infty$ such that for every open set $O \subset
\Rn$,
$$\mu(TO) \leq C\Lambda^{(\infty)}_{n\beta}(O).$$
\end{lemma}
Here $\Lambda^{(\infty)}_d$, $0 < d \leq n$, is $d$-dimensional Hausdorff
capacity, defined (see \cite{A}) for a set $E \subset \Rn$ by
$$\Lambda^{(\infty)}_d(E) = \inf \sum_j r_j^{d},$$
where the infimum is taken over all countable covers $E \subset \bigcup_j
B_j$, with each $B_j$ a ball of radius $r_j$.\\

This paper is concerned with investigating and showing that there
exist similar relationships between Carleson-type measures
and classes of measures generated by capacities.  In addition to
Hausdorff capacity, we consider classes of positive measures on the
half-space generated by Riesz, Bessel, or Besov capacities.

For $0 < \alpha < n$, $1\leq p<n/\alpha$ denote by $h_p^\alpha(\Rn)$
the Riesz potential spaces, defined by the
Riesz potentials, $R^\alpha,$
$$
(R^\alpha g)(x)=c(\alpha,n) \int_{{\bf R}^n} |x-y|^{\alpha-n} g(y)
dy,
$$
where $c(\alpha,n)$ is a certain positive constant (see \cite{M}),
as follows:
$$
h^\alpha_p :=\{f: f=R^\alpha g, g\in L^p\},
$$
with a norm
$$
\|f\|_{h^\alpha_p} :=\|g\|_{L^p}.
$$
Here $L^p$ is the Lebesgue space on $\mathbf {R}^n$.

For any open set $O\subset \mathbf {R}^n,$ the Riesz capacity of
$O$ is defined by (see \cite{AH}),
$$
cap(O;h^\alpha_p):=\inf \{\|f\|_{h^\alpha_p}^p : f\in {\cal S},
f\geq 0, f(x)\geq 1 \;\text{if}\; x\in O\},
$$
where $\cal S$ is the Schwartz class of functions on $\mathbf
{R}^n.$

Let $\Phi$ be a positive function on $(0,\infty),$
$\Phi(0)=0.$
We define the classes of measures $Cap^\Phi(h^\alpha_p)$ as
follows:
$$
\mu\in Cap^\Phi(h^\alpha_p) \quad\text{iff}\quad \mu(TO)\prec
\Phi(cap(O;h^\alpha_p)),$$
uniformly for all open sets $O\subset \mathbf {R}^n$.

Here the
short notation $X\prec Y$ is used for the estimate $X\leq c Y.$ If
$X\prec Y$ and $Y\prec X$ then we write $X\approx Y.$

This definition is the analogue, for measures on the upper
half-space and tents over sets, of a definition given by Maz'ya for
measures and sets in $\Rn$ (see Ch.\ 8 of \cite{M}, also \cite{M2} for
the original idea). Namely, define the class of Maz'ya
measures, $M^\Phi(h^\alpha_p),$ consisting of all positive
measures $\nu$ in $\mathbf{R}^n$ such that
$$\nu(O)\prec \Phi(cap(O;h^\alpha_p)).$$
Note that if $\nu\in
M^\Phi(h^\alpha_p)$ and
$\delta_\epsilon$ denotes the unit mass at some $\epsilon>0$, then
the measure $\mu=\nu\times \delta_\epsilon$ on
$\Rnp$ belongs to $Cap^\Phi(h^\alpha_p)$, since
$$\mu(TO) = \nu(O_\epsilon) \leq \nu(O) \prec \Phi(cap(O;h^\alpha_p)),$$
where $O_\epsilon = \{x \in O: \dist(x, \partial O) > \epsilon\}$.

The main goal of this paper is to give a more geometric characterization
of the classes $Cap^\Phi(h^\alpha_p)$  in terms of
Carleson-type measures $C^\Phi,$ defined as follows:
$$
\mu\in C^\Phi  \quad\text{iff}\quad \mu(TB_r)\prec \Phi(r),
$$
uniformly for all balls $B_r$ in $\mathbf {R}^n$ of radius $r$.
For example, if $\Phi(r)=r^{n\beta},$ $\beta>0,$ then $C^\Phi$ is
the class $C^\beta$ of $\beta-$Carleson measures, as in
(\ref{Cbeta}).

We now state a typical result, proved in the paper. Suppose that
$\Phi$ is a function from $[0,\infty)$ onto itself which is equivalent
to a strictly increasing function $\widetilde{\Phi}$, i.e. $\Phi \approx
\widetilde{\Phi}$, and which satisfies the following conditions:
\begin{equation}
\label{eq00}
\Phi(0)=0, \;
\Phi(cs)\approx \Phi(s),\; c>0,
\end{equation}
(with constants which may depend on $c$), and for some $(\alpha,p),$
$1\leq p<\infty,$ $0<\alpha<n/p$,
\begin{equation}\label{eq01}
\int_0^s \left[ \frac{\Phi(u)}{u^{n-\alpha p}}\right]^{1/p} \frac{du}{u}
\prec\left[
\frac{\Phi(s)}{s^{n-\alpha p}}\right]^{1/p} \quad \mbox{for all} \; s > 0.
\end{equation}
Then
\begin{equation}\label{eq1}
Cap^{\Phi_\alpha} (h^\alpha_p)=C^\Phi\;\text{if}\;
\Phi_\alpha(s)=\Phi(s^{\frac{1}{n-\alpha p}}).
\end{equation}

For example, if $\Phi(s)=s^{n\beta}(1+|\log s|)^\gamma,$  then
$\Phi$ satisfies (\ref{eq00}) and (\ref{eq01}) for $1\leq
p<\infty$, $\beta>1-\alpha p/n$,  $0<\alpha<n/p$, and all real $\gamma$.
Note that if $\gamma \geq 0$ then $\Phi$ is strictly increasing, and
otherwise $\Phi$ is equivalent to a strictly increasing
function $\widetilde{\Phi}$.

In proving (\ref{eq1}), we use the property
\begin{equation}\label{eq2}
cap(B_r; h^\alpha_p)=c r^{n-\alpha p},\; 1\leq p<n/\alpha,
\end{equation}
where $c$ is the capacity of the unit ball. Hence the embedding
\begin{equation}\label{eq3}
Cap^{\Phi_\alpha} (h^\alpha_p)\subset C^\Phi\; 1\leq
p<n/\alpha,\;\Phi_\alpha(s)=\Phi((s/c)^{\frac{1}{n-\alpha p}}),
\end{equation}
is always true.

In order to see the inverse embedding, we prove weak-type
estimates for certain convolution operators of the form
\begin{equation}\label{eq4}
f \mapsto u,\; \; u(x,t) = f\ast \phi_t(x),\; \;
\phi_t(x)=t^{-n}\phi(x/t),
\end{equation}
for certain kernels $\phi$, when the domain is $h^\alpha_p$
and the range is an appropriate Lorentz space in $\mathbf
{R}^{n+1}_+$ with respect to the measure $\mu\in C^\Phi.$ These
weak-type estimates are proved using an analogy with the technique
applied by Adams (see Adams' proof of Theorem 2, Section 1.4.1
of \cite{M}) in proving sharp embeddings of $h^\alpha_p$ into the
Lebesgue space built-up over a positive measure in $\mathbf
{R}^n$.  Moreover, such type of estimates are important for
applications (see \cite{KX}, \cite{X1} and \cite{X2}), so we devote
Section 5  to deriving the corresponding strong-type estimates. They are
related to the problem of sharp embeddings. Because of the analogy between
our classes of measures and Maz'ya measures, we can try to extend the
embedding results (and methods) from
\cite{M} to certain convolution operators of type (\ref{eq4}) and
measures $\mu\in Cap^\Phi(h^\alpha_p)$.

For related embedding theorems in the case of measures on the unit
ball in $\C^n$, see \cite{CO}.

\section{Measures generated by homogeneous capacities}

\subsection{Riesz capacities}

The main result is the following theorem

\begin{theorem}\label{thm1}
Let the function $\Phi$ satisfy conditions (\ref{eq00}) and (\ref{eq01}).
Then
\begin{equation}
\label{eq21}
Cap^{\Phi_\alpha} (h^\alpha_p)=C^\Phi\;\text{if}\;
\Phi_\alpha(s)=\Phi(s^{\frac{1}{n-\alpha p}}).
\end{equation}
\end{theorem}

\begin{proof}
One embedding is given by (\ref{eq3}). To see the other, let
$\mu\in C^\Phi.$ As it has been already explained, we are going to
derive weak-type estimates for the convolution operator $f \mapsto
u$, $u(x,t) = f\ast\phi_t(x),$ where $\phi$ is a positive function such
that (see \cite{KX})
\begin{equation}\label{eq02}
\phi(x)\prec (1+|x|)^{-n},\; \phi\in L^1.
\end{equation}

The proof follows the scheme of the proof of Lemma 3.2 in \cite{KX}.
For
$\lambda>0$, let
\begin{equation}
\label{eq22}
h(\lambda)=\mu(E_\lambda),\quad E_\lambda=\{(x,t) \in \Rnp:
|f\ast\phi_t(x)|>\lambda\} .
\end{equation}
If $f=R^\alpha g,$ $g\in L^p,$ then
$$
f\ast \phi_t(x)=\int_{\mathbf {R}^n} \psi(t,x-z)g(z)dz,
$$
where
$$
\psi(t,z)=c(\alpha,n)\int_{\mathbf {R}^n} |y-z|^{\alpha-n}
\phi_t(y) dy,
$$
and according to Lemma 3.1 of \cite{KX},
\begin{equation}\label{eq23}
\psi(t,x)\prec (t+|x|)^{\alpha-n},\; 0<\alpha<n.
\end{equation}
Therefore
$$
\lambda h(\lambda) \prec \int_{E_\lambda} |f\ast \phi_t(x)|d\mu
\prec\int_{\mathbf {R}^n} |g(z)|\int_{E_\lambda}
(t^2+|x-z|^2)^{\frac{\alpha-n}{2}} d\mu dz,
$$
whence (cf. \cite{KX}, proof of Lemma 3.2)
\begin{equation}\label{eq24}
\lambda h(\lambda) \prec \int_0^\infty r^{\alpha-n-1}
\int_{\mathbf {R}^n} |g(z)|\mu_\lambda (TB(z,r)) dz dr,
\end{equation}
where $\mu_\lambda$ is the restriction of $\mu$ to $E_\lambda.$

Let
$$
T_1(s)=\int_0^s r^{\alpha-n-1} \int_{\mathbf {R}^n}
|g(z)|\mu_\lambda (TB(z,r)) dz dr,
$$
$$
T_2(s)=\int_s^\infty r^{\alpha-n-1} \int_{\mathbf {R}^n}
|g(z)|\mu_\lambda (TB(z,r)) dz dr.
$$

We estimate $T_1(s)$ for arbitrary $s>0$. Using
H\"{o}lder's inequality, we get
$$
T_1(s)\leq \int_0^s r^{\alpha-n-1} \left(\int_{\mathbf {R}^n}
|g(z)|^p \mu_\lambda (TB(z,r)) dz \right)^{1/p}
\left(I(r)\right)^{1 - 1/p}dr,
$$
where $I(r)=\int_{\mathbf {R}^n} \mu_\lambda
(TB(z,r)) dz.$
Since
\begin{equation}\label{eq25}
I(r)=\int_{\mathbf {R}^n} \int_{TB(z,r)\cap E_\lambda} d\mu
dz = \int_{(x,t) \in E_\lambda} \left(\int_{B(x,r-t)}dz\right) d\mu
\prec
r^n h(\lambda),
\end{equation}
and $\mu\in C^\Phi,$ we get
\begin{equation}\label{eq26}
T_1(s)\prec \|g\|_{L^p} [h(\lambda)]^{1 - 1/p} \int_0^s r^{\alpha-n/p}
\Phi^{1/p}(r)
\frac{dr}{r} .
\end{equation}

Analogously, now using the estimate $\mu_\lambda(TB(z,r))\leq
h(\lambda),$ we derive
\begin{equation}\label{eq27}
T_2(s)\prec \|g\|_{L^p}  h(\lambda) s^{\alpha-n/p}\;\;\text{if}\;
\alpha<n/p,\; 1 \leq p < \infty.
\end{equation}
In this way, using also (\ref{eq01}), we have
\begin{equation}\label{eq28}
\lambda h(\lambda)\prec
\|f\|_{h^\alpha_p}\{[h(\lambda)]^{1 - 1/p} [\Phi(s)]^{1/p}
s^{\alpha-n/p}   + h(\lambda) s^{\alpha-n/p} \}.
\end{equation}
Replacing, if necessary,  $\Phi$ by an equivalent strictly increasing
function $\widetilde{\Phi}$ in (\ref{eq28}), we can choose
$s>0$ such that $\Phi(s)=h(\lambda).$ Then
\begin{equation}\label{eq29}
\lambda \left[\Phi^{-1}(h(\lambda))\right]^{n/p-\alpha}\prec
\|f\|_{h^\alpha_p},
\end{equation}
which is the desired weak-type estimate.

Next we prove that (\ref{eq29}) implies $\mu\in
Cap^{\Phi_\alpha}(h^\alpha_p).$ Let $O\subset \mathbf {R}^n$ be an
open set and let $f\in \cal S,$ $f\geq 0,$ $f\geq 1$ on $O.$ If
$(x,t)\in TO,$ then $B(x,t)\subset O,$ hence
$$
f\ast\phi_t(x) =\int_{\mathbf {R}^n} f(z) \phi_t(x-z)dz\geq
\int_{B(x,t)} \phi_t(x-z)dz=\int_{\{z:|z|<1\}} \phi(z)dz.
$$
Therefore,
$$
TO\subset \{(x,t): f\ast\phi_t(x)>d_\phi:=\int_{\{z:|z|<1\}}
\phi(z)dz\}.
$$
Thus
\begin{equation}\label{eq231}
\mu(TO)\prec h(d_\phi).
\end{equation}

On the other hand, the above $f$ can be chosen so that
$\|f\|_{h^\alpha_p}^p\leq 2 cap(O;h^\alpha_p).$ Then (\ref{eq29})
gives
\begin{equation}\label{eq232}
\left[\Phi^{-1}(h(d_\phi))\right]^{n-\alpha p}\prec
cap(O;h^\alpha_p),
\end{equation}
which together with (\ref{eq00}) implies
\begin{equation}\label{eq233}
h(d_\phi)\prec
\Phi\left([cap(O;h^\alpha_p)]^{\frac{1}{n-\alpha p}}\right).
\end{equation}
Note that if we used the equivalent strictly increasing function
$\widetilde{\Phi}$ instead of $\Phi$ in  (\ref{eq29}), we can now go back
to using
$\Phi$ on the right-hand-side.
Combining (\ref{eq231}) and (\ref{eq233}) gives the desired conclusion:
$$
\mu(TO)\prec
\Phi_\alpha(cap(O;h^\alpha_p)),\;\Phi_\alpha(s)=\Phi(s^{\frac{1}{n-\alpha
p}}).
$$
\end{proof}

\subsection{Homogeneous Sobolev capacities}

The homogeneous Sobolev spaces, $w^m_p,$ $1\leq p<\infty,$
$m-$positive integer, are defined as the closure of
$C_0^\infty-$functions on $\mathbf {R}^n$ with respect to the norm
$$
\|f\|_{w^m_p} =\sum_{|\kappa|=m} \|D^\kappa f\|_{L^p}.
$$
As in the previous subsection we can define the homogeneous
Sobolev capacities
$$
cap(O;w^m_p)=\inf\{\|f\|_{w^m_p}^p: f\in {\cal S},\; f\geq 0,\;
f(x)\geq 1\; \text{if}\; x\in O\},
$$
and the classes of measures $Cap^\Phi(w^m_p),$
$$
\mu\in Cap^\Phi(w^m_p) \quad\text{iff}\quad \mu(TO)\prec
\Phi(cap(O;w^m_p)).
$$

Then Theorem \ref{thm1} implies the following
\begin{corollary}\label{cor1}
Let $\Phi$ satisfy the conditions (\ref{eq00}) and (\ref{eq01}).
If $1 < p < n/m$, or $p = 1$ and $m < n$ is even, then
\begin{equation}\label{eq234}
Cap^{\Phi_m}(w^m_p)=C^\Phi,\quad
\Phi_m(s)=\Phi(s^{\frac{1}{n-mp}}).
\end{equation}
\end{corollary}

Indeed, using the Fourier transform ${\cal F}$, defined on functions $f
\in L^1(\Rn)$ by
$${\cal F}f(\xi) = \int_{\Rn} f(x)e^{-ix\cdot \xi}\; dx,$$
one can write the Riesz potential of a function $g \in {\cal S}(\Rn)$ as
$$
R^\alpha g={\cal F}^{-1}(|\xi|^{-\alpha}{\cal F}g),
$$
while the partial derivatives of $f \in {\cal S}(\Rn)$ are given by
$$
\partial_j f = {\cal F}^{-1}(i \xi_j{\cal F}f).
$$
If $m$ is an even integer, and $f = R^m g$, we have
$$g = {\cal F}^{-1}(|\xi|^m{\cal F}f) = {\cal
F}^{-1}((\sum \xi_j^2)^{m/2}{\cal F}f),$$
which leads to the inclusion $w^m_p\subset h^m_p$, $p \geq 1$.

This containment and the opposite one for any integer $m$
and $1 < p < \infty$, follows from the fact that for this range of
$p$ the Riesz transforms
${\cal R}_j$, $1 \leq j \leq n$, defined (for $g \in {\cal S}(\Rn)$) by
$${\cal R}_j g = {\cal F}^{-1}\Big(\frac{i \xi_j}{|\xi|}
{\cal F}g\Big),$$
extend to bounded operators from $L^p$ to $L^p$.
If $f = R^m g$, $m$ an odd integer, then
$$g = {\cal F}^{-1}(|\xi||\xi|^{(m-1)}{\cal F}f) = -i{\cal
F}^{-1}(\sum \frac{i \xi_j}{|\xi|}\xi_j(\sum \xi_j^2)^{(m-1)/2}{\cal
F}f),$$
which gives again $w^m_p\subset h^m_p$ for $1 < p < \infty$.

For the converse, note that for any integer $m$, if $f = R^m g$ with
$g \in L^p(\Rn)$
and $\kappa$ is a multi-index with $|\kappa| = m$, we get
$$D^\kappa f= {\cal F}^{-1}\Big(\frac{(i \xi)^\kappa}{|\xi|^m}
{\cal F}g\Big) = ({\cal R}_1,{\cal R}_2, \ldots,{\cal R}_n)^\kappa g
\in L^p(\Rn), \quad 1 < p < \infty.
$$

Thus if $1<p<n/m,$ then $w^m_p=h^m_p$ and (\ref{eq234})
follows from Theorem \ref{thm1}. If $p=1$ and $m$ is even then
$w^m_1\subset h^m_1,$ hence
$$
C^\Phi=Cap^{\Phi_m}(h_1^m)\subset Cap^{\Phi_m}(w_1^m).
$$
Finally, since
$$
cap(B_r;w_1^m)=c r^{n-m},
$$
we have
$$
Cap^{\Phi_m}(w_1^m)\subset C^\Phi.
$$

\subsection{Homogeneous Besov capacities}

We can define the homogeneous Besov spaces, $b^\alpha_{p,q},$
$\alpha>0,$ $1\leq p<\infty,$ $0<q\leq \infty,$ by interpolation:
\begin{equation}\label{eq235}
b^\alpha_{p,q}=(L^p, w^m_p)_{\alpha/m,q},\;0<\alpha<m,
\end{equation}
where $(\cdot, \cdot)_{\sigma,q}$ stands for the real
interpolation method (see, for example, \cite{BL}).

For any open set $O\subset {\bf R}^n,$ the homogeneous Besov
capacity of $O$ is defined by
$$
cap(O;b^\alpha_{p,q})=\inf\{\|f\|_{b^\alpha_{p,q}}^p: f\in {\cal
S},\; f\geq 0,\; f(x)\geq 1\; \text{if}\; x\in O\},
$$
and the classes of measures $C^\Phi(b^\alpha_{p,q})$ by
$$
\mu\in C^\Phi(b^\alpha_{p,q})\;\text{iff}\; \mu(TO)\prec
\Phi(cap(O;b^\alpha_{p,q})),\;\text{uniformly for all open sets}\;
O\subset \mathbf {R}^n.
$$

\begin{theorem}\label{thm2}
Let the function $\Phi$ satisfy conditions (\ref{eq00}) and (\ref{eq01}).
Then
\begin{equation}\label{eq236}
Cap^{\Phi_\alpha}(b^\alpha_{p,q})=C^\Phi\;\text{if}\;1<p<n/\alpha,\;
0<q\leq\infty,\; \Phi_\alpha(s)=\Phi(s^{\frac{1}{n-\alpha p}}).
\end{equation}
\end{theorem}

\begin{proof}
Since
$$
cap(B_r;b^\alpha_{p,q})=c r^{n-\alpha p},\; 1\leq p<n/\alpha,
$$
where $c=cap(B_1;b^\alpha_{p,q}),$ we have the embedding
$$
Cap^{\Phi_\alpha}(b^\alpha_{p,q})\subset C^\Phi\;\text{if}\; 1\leq
p<n/\alpha\;\text{and}\; \Phi_\alpha(s)=\Phi((s/c)^{\frac{1}{n-\alpha
p}}),
$$
for any $\Phi>0.$

To see the inverse, we use the weak-type estimate (\ref{eq29}) for $\mu\in
C^\Phi$ and real interpolation for fixed $p>1$ and $\Phi.$ Since
(\cite{BL}, Theorem 6.3.1),
$$
(h^{\alpha_1}_p, h^{\alpha_2}_p)_{\theta,q}=b^\alpha_{p,q},\;
1<p<\infty,\; 0<q\leq \infty,\; \alpha=(1-\theta)\alpha_1+\theta
\alpha_2,\;0<\theta<1,
$$
we derive from (\ref{eq29}),
\begin{equation}\label{eq237}
\lambda \left[\Phi^{-1}(h(\lambda))\right]^{n/p-\alpha} \prec
\|f\|_{b^{\alpha}_{p,q}}.
\end{equation}
Note that if we had to use the equivalent strictly increasing function
$\widetilde{\Phi}$ instead of $\Phi$ in  (\ref{eq29}), we would now have
$\widetilde{\Phi}^{-1}$ in (\ref{eq237}).
As at the end of the proof of Theorem \ref{thm1}, we conclude that
(\ref{eq237}) implies $\mu\in Cap^{\Phi_\alpha}(b^\alpha_{p,q}).$
\end{proof}

\begin{theorem}[Case $b^\alpha_{1,q}$]\label{thm3}
Let the function $\Phi$ satisfy conditions (\ref{eq00}) and
\begin{equation}\label{eq011}
\int_0^s \left(\frac{\Phi(u)}{u^n}\right)^{1/r} \frac{du}{u}\prec
\left(\frac{\Phi(s)}{s^n}\right)^{1/r},\;1/r=1-\alpha/n,\;0<\alpha<n.
\end{equation}
Then
\begin{equation}\label{eq2361}
Cap^{\Phi_\alpha}(b^\alpha_{1,q})=C^\Phi, \;
\Phi_\alpha(s)=\Phi(s^{\frac{1}{n-\alpha }}).
\end{equation}
\end{theorem}
Note that (\ref{eq011}) implies (\ref{eq01}) if $p=1.$

\begin{proof}
If $\mu\in C^\Phi$, $f\in L^r,$ and $\Phi$ satisfies
(\ref{eq011}), then the same proof as that of Theorem \ref{thm1}
(formally taking $\alpha=0$), and again replacing $\Phi$ by an
equivalent strictly increasing function if necessary, shows that
$$
\lambda \left(\Phi^{-1}(h(\lambda))\right)^{n/r}\prec
\|f\|_{L^r},\; 1\leq r<\infty.
$$
We can interpolate this inequality, hence
$$
\lambda \left(\Phi^{-1}(h(\lambda))\right)^{n/r}\prec
\|f\|_{L^{r,q}},\; 1< r<\infty, 0<q\leq \infty,
$$
where $L^{r,q}$ is the Lorentz space (see \cite{BL}).

Since we have the embedding
$$
b^\alpha_{1,q}\subset L^{r,q},\; 1/r=1-\alpha/n, 0<\alpha<n,
$$
we get
\begin{equation}\label{eq2371}
\lambda \left(\Phi^{-1}(h(\lambda))\right)^{n/r}\prec
\|f\|_{b^{\alpha}_{1,q}},
\end{equation}
i.e.\ the estimate (\ref{eq237}) for $p=1.$ As before, we conclude
from (\ref{eq2371}) that $\mu \in
Cap^{\Phi_\alpha}(b^{\alpha}_{1,q}).$
\end{proof}

\begin{remark}
\label{remark21}
The same proofs as those of Theorems \ref{thm1} and \ref{thm2}
show the following embeddings for a larger classes of functions
$\Phi.$ Namely, let $\Phi$ satisfy (\ref{eq00}) and let
\begin{equation}\label{eq238}
F(s):=s\left(\int_0^s
\left(\frac{\Phi_\alpha(u)}{u}\right)^{1/p}\frac{du}{u}\right)^p,\;
1\leq p<n/\alpha,\; \Phi_\alpha(s)=\Phi(s^{\frac{1}{n-\alpha p}}).
\end{equation}
Then
\begin{equation}\label{eq239}
C^\Phi\subset Cap^F(h^\alpha_p)\; 1\leq p<n/\alpha
\end{equation}
and
\begin{equation}\label{eq240}
C^\Phi\subset Cap^F(b^{\alpha}_{p,q})\; 1<p<n/\alpha.
\end{equation}
\end{remark}

For example, if $\Phi(u)=u^{n-\alpha p} (1+|\log u|)^\gamma,$
$\gamma<-p,$ then $ F(s)\approx s(1+|\log s|)^{\gamma+p}.$

\section{Measures generated by inhomogeneous capacities}

\subsection{Bessel capacities}

We first recall the definition of the Bessel potential spaces,
$H^\alpha_p,$ $0<\alpha<n,$ $1\leq p<\infty$ (see \cite{AH},
\cite{M}). We say that $f\in H^\alpha_p$ iff $f=G_\alpha\ast g,$
$g\in L^p,$ and the norm is given by
$\|f\|_{H^\alpha_p}=\|g\|_{L^p},$ where $G_\alpha$ is the Bessel
kernel
$$
G_\alpha(x)=|x|^{(\alpha-n)/2} K_{\frac{n-\alpha}{2}}(|x|),
$$
and $K$ is the modified Bessel function of third kind:
$$
K_{\frac{n-\alpha}{2}}(|x|)=c |x|^{-1/2} e^{-|x|} \int_0^\infty
e^{-u} u^{\frac{n-\alpha-1}{2}}
\left(1+\frac{u}{2|x|}\right)^{\frac{n-\alpha-1}{2}} du.
$$
In particular, we have the global estimate
\begin{equation}\label{eq31}
G_\alpha(x)\prec e^{-|x|} \left(1+|x|^{\alpha-n}\right).
\end{equation}

For any open set $O\subset \mathbf {R}^n,$ we define its Bessel
capacity by
$$
cap(O;H^\alpha_{p})=\inf\{\|f\|_{H^\alpha_{p}}^p: f\in {\cal S},\;
f\geq 0,\; f(x)\geq 1\; \text{if}\; x\in O\}.
$$
For example (see \cite{AH}),
\begin{equation}\label{eq33}
cap(B_r;H^\alpha_p)\approx r^{n-\alpha p} \;\text{if}\; 0<r<1,\;
0<\alpha<n/p,\; 1<p<\infty
\end{equation}
and
\begin{equation}\label{eq34}
cap(B_r;H^\alpha_p)\approx r^{n} \;\text{if}\; r>1,\;
0<\alpha<n/p,\; 1<p<\infty.
\end{equation}

The classes of measures $Cap^\Phi(H^{\alpha}_p)$ are defined as
follows:
$$
\mu\in Cap^\Phi(H^{\alpha}_p)\;\text{iff}\;
\mu(TO)\prec\Phi(cap(O;H^{\alpha}_{p})),
$$
uniformly for all open sets $O\subset \mathbf{R}^n.$

Our main result in the nonhomogeneous context is the following theorem.
\begin{theorem}\label{thm31}
Let $\Psi$ be a function from $[0,\infty)$ onto itself which is
equivalent to a strictly increasing function and satisfies  condition
(\ref{eq00}), as well as
\begin{equation}\label{eq012}
\int_0^s \left(\frac{\Psi(u)}{u}\right)^{1/p}
\frac{du}{u}\prec\left(\frac{\Psi(s)}{s}\right)^{1/p}
\end{equation}
for every $s > 0$.  Set
\begin{equation}\label{eq36}
\Psi^\alpha(s)=\left\{ \begin{array}{ll}
                                                        \Psi(s^{n-\alpha p}) & \mbox{ if $0<s<1$}, \\
                                                        \Psi(s^{n}) & \mbox{ if $s\geq1$}.\\
                     \end{array}\right.
\end{equation}
Then
\begin{equation}
Cap^\Psi(H^\alpha_p)= C^{\Psi^\alpha},\; 1< p<n/\alpha.
\end{equation}
\end{theorem}

\begin{proof}
Using (\ref{eq33}), (\ref{eq34}), we obtain the embedding
\begin{equation}\label{eq35}
Cap^\Psi(H^\alpha_p)\subset C^{\Psi^\alpha},\; 1\leq p<n/\alpha.
\end{equation}
Conversely, assume $\mu\in C^{\Psi^\alpha}$. As in the proof of Theorem
\ref{thm1}, we can write
\begin{equation}\label{eq367}
\lambda h(\lambda) \prec \int_{E_\lambda} |f\ast \phi_t(x)|d\mu
\prec\int_{\Rn} \int_{\Rn}|g(z)| G_\alpha\ast
\phi_t(x-z) d\mu_\lambda dz.
\end{equation}

We need the estimate
\begin{equation}\label{eq37}
G_\alpha\ast\phi_t(x)\prec \min\{(t+|x|)^{\alpha-n},
(t+|x|)^{-n}\},
\end{equation}
if $\phi$ is a non-negative smooth function with compact support
in the unit ball, and $d_\phi:=\int_{\mathbf {R}^n} \phi(x)dx>0.$

The first estimate is a consequence of $G_\alpha(x)\prec
|x|^{\alpha-n}$ and (\ref{eq23}).
To prove the second estimate, we use the properties
$G_\alpha\in L^1$ and
$\Phi_t\prec t^{-n}.$ This gives $G_\alpha\ast\phi_t(x)\prec t^{-n}$,
or
\begin{equation}\label{eq38}
G_\alpha\ast\phi_t(x)\prec  (t+|x|)^{-n}\;\text{if}\; |x|<2t.
\end{equation}

If $|x|>2t$ and $y\in \;\text{support of}\; \phi_t$ (hence
$|y|\leq t$), then $|x|+t\prec  |x|/2 \leq |x-y|.$ Therefore
(\ref{eq31}) implies
$$
G_\alpha\ast\phi_t(x)\prec C_N \int_{\mathbf{R}^n} |x-y|^{-N}
(1+|x-y|^{\alpha-n}) \phi_t(y) dy,
$$
or
\begin{equation}\label{eq39}
G_\alpha\ast\phi_t(x)\prec  (t+|x|)^{-n}\;\text{if}\; |x|>2t.
\end{equation}

In estimating $h(\lambda)$, we consider two cases.  We assume that
$\Psi$ is strictly increasing, otherwise we replace it by a strictly
increasing function $\widetilde{\Psi}$ with $\Psi \approx
\widetilde{\Psi}$.

Case 1: $h(\lambda)<\Psi(1)$.
In this case we apply the first estimate in (\ref{eq37}). Thus we
get (\ref{eq24}) so we can argue as in the proof of Theorem \ref{thm1}.
We take $s < 1$ so that in estimate (\ref{eq26}) for $T_1(s)$
we can use the fact that $\mu_\lambda(TB(x,r)) \leq  \Psi(r^{n-\alpha
p})$. With the change of variable $r^{n-\alpha p} \mapsto r$ (note
$\alpha<n/p$), and using (\ref{eq012}), we get
\begin{equation}\label{eq310}
\lambda h(\lambda)\prec \|f\|_{H^\alpha_p}\{
 s^{-1/p} [\Psi(s)]^{1/p} [h(\lambda)]^{1 - 1/p} + s^{-1/p} h(\lambda)\}.
\end{equation}
Since
$h(\lambda)<\Psi(1)$, and we are assuming $\Psi$ is strictly increasing,
we can choose
$s < 1$ such that
$\Psi(s)=h(\lambda).$ In this case we get
\begin{equation}\label{eq311}
\lambda \left[\Psi^{-1}(h(\lambda))\right]^{1/p}\prec
\|f\|_{H^\alpha_p}.
\end{equation}

Case 2: $h(\lambda)\geq \Psi(1)$.  Here we can assume $s > 1$.
We start with (\ref{eq367}) and write
\begin{equation}\label{eq312}
\lambda h(\lambda)\prec I + II + III,
\end{equation}
where
$$
I =\int_{\Rn} \int_{\{(x,t): |x-z|+t<1/2\}} |g(z)| G_\alpha\ast
\phi_t(x-z) d\mu_\lambda dz,
$$
$$
II =\int_{\Rn} \int_{\{(x,t): 1/2<|x-z|+t<s/2\}} |g(z)| G_\alpha\ast
\phi_t(x-z) d\mu_\lambda dz,
$$
$$
III =\int_{\Rn} \int_{\{(x,t): |x-z|+t>s/2\}} |g(z)| G_\alpha\ast
\phi_t(x-z) d\mu_\lambda dz,\; 1\prec s.
$$

The first integral can be estimated in the same way as $T_1(s_0)$
for some constant $s_0$, so that (\ref{eq26}) simplifies to
\begin{equation}\label{eq313}
I\prec \|f\|_{H^\alpha_p}  [ h(\lambda)]^{1-1/p}.
\end{equation}
To estimate $II,$ we use the second bound in (\ref{eq37}), whence
$$
II \prec \int_{\Rn}\int_{\{(x,t): 1/2<\rho(x-z,t)<s/2\}} |g(z)|
[\rho(x-z,t)]^{-n} d\mu_\lambda dz,
$$
where $\rho(x,t):=|x|+t.$ Since
$$
\rho^{-n}\prec \int_\rho^s r^{-n-1}dr\;\text{if}\; \rho<s/2,
$$
now we have, instead of (\ref{eq24}),
$$
II\prec \int_{1/2}^s r^{-n-1} \int_{\mathbf {R}^n}
\mu_\lambda(TB(z,r))dzdr.
$$
Arguing as before and using $\mu(TB(z,r))\prec \Psi^\alpha(r)
\prec \Psi(r^n)$ if
$r \geq 1/2$ (by definition (\ref{eq36}) and property (\ref{eq00})
for $\Psi$),  we
get
$$
II\prec \|f\|_{H^\alpha_p} \int_0^s r^{-n/p-1}\Psi(r^n)^{1/p} dr [
h(\lambda)]^{1-1/p}.
$$
Again changing variables $r^n\mapsto r$, and using
(\ref{eq012}), we obtain
\begin{equation}\label{eq314}
II\prec \|f\|_{H^\alpha_p} s^{-1/p} [\Psi(s)]^{1/p} [
h(\lambda)]^{1-1/p}.
\end{equation}

Analogously,
$$
III\prec \int_{s/2}^\infty r^{-n-1} \int_{\mathbf {R}^n}
\mu_\lambda(TB(z,r))dzdr,
$$
and using $\mu_\lambda(TB(z,r))\leq h(\lambda)$ we get, as before
(also changing  variables $r^n \mapsto r$),
\begin{equation}\label{eq315}
III\prec \|f\|_{H^\alpha_p} s^{-1/p}  h(\lambda).
\end{equation}
Unifying (\ref{eq312}), (\ref{eq313}), (\ref{eq314}) and
(\ref{eq315}), we get (\ref{eq310}) for $s \geq 1$ as well. Since
now $h(\lambda)>\Psi(1)$, and again assuming we've replaced $\Psi$, if
necessary, by an equivalent strictly increasing function, we can choose
$s \geq 1$ to solve the equation $\Psi(s)=h(\lambda)$, and therefore
(\ref{eq311}) holds in this case as well.

In order to prove that (\ref{eq311}) implies $\mu\in
Cap^\Psi(H^\alpha_p),$ we start with (\ref{eq231}). Analogously to
(\ref{eq232}), we derive from (\ref{eq311}) that
$$
\Psi^{-1}(h(d_\phi))\prec cap(O;H^\alpha_p).
$$
Together with (\ref{eq231}) this means that $\mu\in
Cap^\Psi(H^\alpha_p).$
\end{proof}

\subsection{Inhomogeneous Sobolev capacities}

If $\alpha=m$ is integer and $1<p<\infty$, then $H^m_p=W^m_p$ - the
Sobolev space with norm
$$
\|f\|_{W^m_p} =\sum_{|\kappa|\leq m} \|D^\kappa f\|_{L^p}.
$$
As in the homogeneous case, this can be seen via the Fourier
transform, since the Bessel potential  $f=G_m\ast g$ can be
written as
\begin{equation}
\label{Bessel-FT}
f = {\cal F}^{-1} ((1 + |\xi|^2)^{-m/2} {\cal F}g),
\end{equation}
and the operators defined by the Fourier multipliers
$\frac{\xi^\kappa}{(1 + |\xi|^2)^{m/2}}$, $|\kappa| \leq m$, are
bounded on $L^p(\Rn)$ for $1<p<\infty$ (see also \cite{S}, Ch.\ V,
Theorem 3.3).
When $p = 1$ and $m$ is even, (\ref{Bessel-FT}) shows
\begin{equation}\label{eq32}
W^m_1\subset H^m_1,
\end{equation}
but this inclusion fails for $m$ odd, and equality does not hold
(see \cite{S}, Ch.\ V, Section 6.6).

By definition, for any open set $O\subset \mathbf {R}^n,$
$$
cap(O;W^m_{p})=\inf\{\|f\|_{W^m_{p}}^p: f\in {\cal S},\; f\geq
0,\; f(x)\geq 1\; \text{if}\; x\in O\},
$$
and
$$
\mu\in Cap_\Phi(W^m_p)\;\text{iff}\;
\mu(TO)\prec\Phi(cap(O;W^m_{p})),
$$
uniformly for all open sets $O\subset \Rn$.

Theorem \ref{thm31} implies the following corollary.

\begin{corollary}\label{cor31}
Let $m$ be an integer less than $n$.
With $\Psi$ and $\Psi^m$ as in Theorem \ref{thm31}, we have
\begin{equation}\label{eq318}
Cap^\Psi(W^m_p)=
C^{\Psi^m}
\end{equation}
for $1< p<n/m$ or $p = 1$ and $m$ even.
\end{corollary}

\begin{proof}
If $1<p<n/m$ then $W^m_p=H^m_p$ and (\ref{eq318}) follows from
Theorem~\ref{thm31} . If $p=1$ and $m$ is even then the inclusion
(\ref{eq32}) implies
\begin{equation}\label{eq319}
C^{\Psi^m}=Cap^\Psi(H^m_1)\subset Cap^\Psi(W^m_1).
\end{equation}
On the other hand, for $m<n$ we have
$$
cap(B_r;W^m_1)\prec
r^{n-m} \;\text{if}\; 0<r<1
$$
and
$$
cap(B_r;W^m_1)\prec r^{n} \;\text{if}\; r\geq1.
$$
Therefore,
\begin{equation}\label{eq320}
Cap^\Psi(W^m_1)\subset
C^{\Psi^m}, \; 1\leq m<n.
\end{equation}
Combining (\ref{eq319}) and (\ref{eq320}) gives
(\ref{eq318}).
\end{proof}

\subsection{Inhomogeneous Besov capacities}

We can define the inhomogeneous Besov spaces, $B^\alpha_{p,q},$
$\alpha>0,$ $1\leq p<\infty,$ $0<q\leq \infty,$ by interpolation:
\begin{equation}\label{eq321}
B^\alpha_{p,q}:=(L^p, W^m_p)_{\alpha/m,q},\;0<\alpha<m.
\end{equation}
We need the following formula (see \cite{BL})
\begin{equation}\label{eq322}
B^\alpha_{p,q}=(H^{\alpha_1}_p,
H^{\alpha_2}_p)_{\theta,q},\;\alpha=(1-\theta)\alpha_1+\theta\alpha_2,\;
1<p<\infty,\; 0<q\leq \infty.
\end{equation}

For any open set $O\subset \Rn$ the inhomogeneous
Besov capacity of $O$ is defined by
$$
cap(O;B^\alpha_{p,q})=\inf\{\|f\|_{B^\alpha_{p,q}}^p: f\in {\cal
S},\; f\geq 0,\; f(x)\geq 1\; \text{if}\; x\in O\}
$$
and the classes of measures $C^\Phi(B^\alpha_{p,q})$ as follows:
$$
\mu\in C^\Phi(B^\alpha_{p,q})\quad\text{iff}\quad
\mu(TO)\prec\Phi(cap(O;B^\alpha_{p,q})),
$$
uniformly for all open sets $O\subset \Rn$.

\begin{theorem}\label{thm32}
Let $\Omega$ be a function from $[0,\infty)$ onto itself, which is
equivalent to a strictly increasing function, and
satisfies  condition (\ref{eq00}).  For $0<\alpha<n/p$ and
$1<p<\infty$, assume $\Omega$ satisfies condition
(\ref{eq01}) whenever $0<s<1$, while when $s \geq 1$ it satisfies
\begin{equation}\label{eq013}
\int_1^s
\left(\frac{\Omega(u)}{u^n}\right)^{1/p}\frac{du}{u}\prec
\left(\frac{\Omega(s)}{s^n}\right)^{1/p}\quad
\text{and}\quad s^n\prec \Omega(s).
\end{equation}
Then for $0<q\leq\infty$,
\begin{equation}\label{eq323}
Cap^{\Omega_\alpha}(B^\alpha_{p,q})=C^\Omega,
\end{equation}
where
\begin{equation}
\Omega_\alpha(s)=\left\{ \begin{array}{ll}
                                            \Omega(s^{\frac{1}{n-\alpha p}}) & \mbox{ if $0<s<1$}, \\
                                            \Omega(s^{\frac{1}{n}}) & \mbox{ if $s\geq1$}.\\
                     \end{array}\right.
\end{equation}
\end{theorem}

For example, with $\beta_0>1-\alpha p/n$ and  $\beta_1>1$,
we can take
\begin{equation}
\Omega(s)=\left\{ \begin{array}{ll}
                                            s^{n\beta_0}(1+|\log s|)^{\gamma_0} & \mbox{ if $0<s<1$}, \\
                                            s^{n\beta_1}(1+|\log s|)^{\gamma_1} & \mbox{ if $s\geq1$}.\\
                     \end{array}\right.
\end{equation}

\begin{proof}
Let $\Psi(s)=\Omega_\alpha(s).$ Then $\Psi^\alpha(s)=\Omega(s),$
where $\Psi^\alpha$ is defined by (\ref{eq36}).  Moreover, $\Psi$
satisfies (\ref{eq012}) and (\ref{eq00}), and if $\Omega \approx
\widetilde{\Omega}$ for some strictly increasing function
$\widetilde{\Omega}$, then
$\Psi$ is equivalent to the strictly increasing function
$\widetilde{\Omega}_\alpha$. Therefore, with
$\mu\in C^\Omega$, we have the estimate (\ref{eq311}). This estimate
can be interpolated for fixed $\Omega$ and $p$. Using
(\ref{eq322}), we derive from (\ref{eq311})
\begin{equation}
\label{eq324}
\lambda\left[\Psi^{-1}(h(\lambda))\right]^{1/p}\prec
\|f\|_{B^\alpha_{p,q}}.
\end{equation}

From (\ref{eq324}), it follows as before that $\mu\in
Cap^\Psi(B^\alpha_{p,q}),$ i.e.
$$
C^\Omega\subset Cap^{\Omega_\alpha}(B^\alpha_{p,q}).
$$
To see the inverse inclusion, we notice that for $0<\alpha<n/p,$
$cap(B_r;B^\alpha_{p,q})\prec r^{n-\alpha p}$ if $0<r<1$ and
$cap(B_r;B^\alpha_{p,q})\prec r^{n}$ if $r>1$. Hence
$$
Cap^{\Omega_\alpha}(B^\alpha_{p,q})\subset C^\Omega.
$$
\end{proof}

\section{Relation with Hausdorff capacities}

For any open set $O\subset{\bf R}^n$ and any positive increasing
function $w$ on $(0,\infty),$ define the ($w$-)Hausdorff capacity of $O$
by
$$
\Lambda^\infty_w(O)=\inf \sum w(r_j),
$$
where the infimum is taken over all coverings of $O$ by countable
unions of balls of radii $r_j$, $O\subset \cup B_{r_j}.$ In
particular,
\begin{equation}\label{eq41}
\Lambda^\infty_w(B_r)\leq w(r).
\end{equation}

If $w(r)=r^d,$ $d>0,$ this is the $d$-dimensional Hausdorff capacity
(or Hausdorff content) as defined by Adams (see \cite{A}), and we
write $\Lambda^\infty_d$ instead of $\Lambda^\infty_{r^d}.$
Since the set function $O\mapsto cap(O;h^\alpha_p)$ is countably
subadditive (see \cite{AH}, p. 26), we see that
\begin{equation}\label{eq42}
cap(O;h^\alpha_p)\prec \Lambda^\infty_{n-\alpha
p}(O),\;1<p<n/\alpha.
\end{equation}

The inverse inequality can not be true uniformly for all open sets
$O$ (see \cite{AH}, p. 148). Using the Hausdorff capacities, we
can define classes of positive measures
$Cap^\Phi(\Lambda_w^\infty)$ in $\mathbf {R}^{n+1}_+$ as follows:
$$
\mu\in Cap^\Phi(\Lambda_w^\infty)\;\text{iff}\; \mu(TO)\prec
\Phi(\Lambda^\infty_w(O)),
$$
uniformly for all open sets $O\subset \mathbf {R}^n$.

\begin{theorem}\label{thm41}
Let the function $\Phi$ be equivalent to a strictly increasing
function and satisfy
conditions (\ref{eq00}) and (\ref{eq01}). Then
\begin{equation}\label{eq43}
Cap^{\Phi_\alpha}(\Lambda^\infty_{n-\alpha p})=Cap^{\Phi_\alpha}
(h^\alpha_p)=C^\Phi\;\text{if}\;
\Phi_\alpha(s)=\Phi(s^{\frac{1}{n-\alpha p}}).
\end{equation}
\end{theorem}

\begin{proof}
From (\ref{eq42}) we derive the embedding
$$
Cap^{\Phi_\alpha} (h^\alpha_p)\subset Cap^{\Phi_\alpha}
(\Lambda^\infty_{n-\alpha p}),\; 1<p<n/\alpha,
$$
and using (\ref{eq41}) we see that
$$
Cap^{\Phi_\alpha} (\Lambda^\infty_{n-\alpha p})\subset
C^\Phi,\;\Phi_\alpha(s)=\Phi(s^{\frac{1}{n-\alpha p}}).
$$
It remains to apply Theorem \ref{thm1}.
\end{proof}

Note that the equality $Cap^{\Phi_\alpha}(\Lambda^\infty_{n-\alpha
p}) = C^\Phi$ for $\Phi(r) = r^{n - \alpha p}$ is just
Lemma~\ref{lemma-Cbeta} with $\beta = 1 - \alpha p /n$.

Analogously, for Bessel capacities we have
\begin{equation}\label{eq44}
cap(O;H^\alpha_p)\leq \Lambda^\infty_{w_{\alpha}}(O),\;1<p<n/\alpha,
\end{equation}
where $w_\alpha(r)=r^{n-\alpha p}$ if $0<r<1,$ and
$w_\alpha(r)=r^n$ if $r>1.$

\begin{theorem}\label{thm42}
Let the function $\Psi$ satisfy conditions (\ref{eq00}) and
(\ref{eq013}). Then
\begin{equation}\label{eq45}
Cap^{\Psi}(\Lambda^\infty_{w_\alpha})=Cap^{\Psi}
(H^\alpha_p)=C^{\Psi^\alpha}\;\text{if}\; \Psi^\alpha=\Psi \circ
w_\alpha.
\end{equation}
\end{theorem}

\begin{proof}
From (\ref{eq44}) it follows that
$$
Cap^{\Psi} (H^\alpha_p)\subset
Cap^{\Psi}(\Lambda^\infty_{w_\alpha}),
$$
and using (\ref{eq41}) we derive
$$
Cap^{\Psi}(\Lambda^\infty_{w_\alpha})\subset C^{\Psi^\alpha}.
$$
It remains to apply Theorem \ref{thm31}.
\end{proof}

\section{Related convolution operators}

Here we prove strong type estimates for the convolution operators
(\ref{eq4}). Let $\phi$ be positive function on $\mathbf {R}^n$,
satisfying the following condition (see \cite{S2}, Chapter II, Section
2.4):
\begin{equation}\label{eq014}
\phi\;\text{ has a non-increasing radial majorant
that is integrable and bounded.}
\end{equation}

For any positive strictly increasing function $\Phi$ on $(0,\infty)$ and
any positive measure $\mu$ on $\mathbf {R}^{n+1}_+$, let
$\Lambda_\mu^p(\Phi^{-1}),$ $1\leq p<\infty,$ denote the Lorentz
space on ${\bf R}^{n+1}_+$, consisting of all measurable functions
$F(x,t)$ such that
$$
\|F\|_{\Lambda^p_\mu(\Phi^{-1})} =\left(\int_0^\infty
\Phi^{-1}(h(\lambda))d\lambda^p\right)^{1/p}<\infty,
$$
where $h(\lambda):=\mu\{(x,t): |F(x,t)|>\lambda\}.$
If $\Phi$ is not strictly increasing but only equivalent to a strictly
increasing function $\widetilde{\Phi}$, we replace $\Phi^{-1}$ by
$\widetilde{\Phi}^{-1}$ in the above definition (where the size of the
norm may depend on the choice $\widetilde{\Phi}$),
but for the sake of simplicity we keep the same notation.

\begin{theorem}\label{thm51}
Let $\phi$ satisfy (\ref{eq014}) and let $\Phi$ be as above and satisfy
(\ref{eq00}). If $\mu\in Cap^\Phi(h^\alpha_p),$ then
\begin{equation}\label{eq46}
\|f\ast\phi_t\|_{\Lambda^p_\mu(\Phi^{-1})}\prec
\|f\|_{h^\alpha_p}, \;1<p<n/\alpha.
\end{equation}
Conversely, the estimate (\ref{eq46}) implies $\mu\in
Cap^\Phi(h^\alpha_p).$
\end{theorem}

\begin{proof}
We start with (\ref{eq22}) and use the relation (see \cite{S2})
\begin{equation}\label{eq47}
h(\lambda)\prec \mu(T\{x: Mf(x)>c\lambda\}),
\end{equation}
where $M$ is the Hardy-Littlewood maximal function and $c$ is a positive
constant depending on $\phi.$

If $\mu\in Cap^\Phi(h^\alpha_p),$ then (\ref{eq47}) implies
\begin{equation}\label{eq48}
\Phi^{-1}(h(\lambda))\prec cap(\{x: Mf(x)>c\lambda\};h^\alpha_p).
\end{equation}
Hence
\begin{equation}\label{eq49}
\int_0^\infty \Phi^{-1}(h(\lambda))d\lambda^p \prec \int_0^\infty
\ cap(\{x: Mf(x)>c\lambda\};h^\alpha_p) d\lambda^p.
\end{equation}
It remains to apply Dahlberg's estimate (see  \cite{M}).
\end{proof}

As a consequence of Theorems \ref{thm1} and \ref{thm51}, we get
(see also \cite{KX}, where the case $\Phi(s)=s^{n\beta},$
$\beta>1-\alpha n/p$ is covered)
\begin{corollary}\label{cor51}
Let $\phi$ satisfy (\ref{eq014}) and let $\Phi$ satisfy
(\ref{eq00}) and (\ref{eq01}). If $\mu\in C^{\Phi_\alpha},$
$\Phi_\alpha(s)=\Phi(s^{n-\alpha p}),$ $1<p<n/\alpha,$ then
\begin{equation}\label{eq410}
\|f\ast\phi_t\|_{\Lambda^p_\mu(\Phi^{-1})}\prec
\|f\|_{h^\alpha_p}.
\end{equation}
\end{corollary}

\begin{theorem}\label{thm52}
Let $\phi$ satisfy (\ref{eq014}) and let $\Phi$ satisfy
(\ref{eq00}). If $\mu\in Cap^\Phi(b^\alpha_{p,q}),$ then
\begin{equation}\label{eq411}
\|f\ast\phi_t\|_{\Lambda^a_\mu((\Phi^{-1})^b)}\prec
\|f\|_{b^\alpha_{p,q}}, \;1\leq p<n/\alpha,\; 1<q\leq \infty,
\end{equation}
where
\begin{equation}\label{eq412}
a=\max\{p,q\},\; b=\max\{1,q/p\}.
\end{equation}
Conversely, the estimate (\ref{eq411}) implies $\mu\in
Cap^\Phi(b^\alpha_{p,q}).$
\end{theorem}

\begin{proof}
For $\mu\in Cap^\Phi(b^\alpha_{p,q})$ we have, analogously to
(\ref{eq49}),
$$
\int_0^\infty[\Phi^{-1}(h(\lambda))]^b d\lambda^a \prec
\int_0^\infty \ [cap(\{x: Mf(x)>c\lambda\};b^\alpha_{p,q})]^b
d\lambda^a.
$$
Applying Corollary 1 of \cite{AX} and Lemma 4.1 of \cite{X2}, we get
(\ref{eq411}).
\end{proof}

As a consequence of Theorems \ref{thm2} and \ref{thm52}, we get
(see also \cite{KX})
\begin{corollary}\label{cor52}
Let $\phi$ satisfy (\ref{eq014}) and let $\Phi$ satisfy
(\ref{eq00}) and (\ref{eq01}). If $\mu\in C^{\Phi_\alpha},$
$\Phi_\alpha(s)=\Phi(s^{n-\alpha p}),$ $1<p<n/\alpha,$ then
\begin{equation}\label{eq413}
\|f\ast\phi_t\|_{\Lambda^a_\mu((\Phi^{-1})^b)}\prec
\|f\|_{b^\alpha_{p,q}},
\end{equation}
where $a$, $b$ are defined by (\ref{eq412}).
\end{corollary}

Using Theorems \ref{thm51} and \ref{thm52}, and Remark
\ref{remark21}, we have the following estimates for a larger
classes of functions $\Phi.$
\begin{corollary}\label{cor525}
Let $\phi$ satisfy (\ref{eq014}) and let $\Phi$ satisfy
(\ref{eq00}). If $\mu\in C^{\Phi},$  then
\begin{equation}\label{eq4131}
\|f\ast\phi_t\|_{\Lambda^p_\mu(F^{-1})}\prec
\|f\|_{h^\alpha_{p}},\; 1<p<n/\alpha
\end{equation}
and
\begin{equation}\label{eq4132}
\|f\ast\phi_t\|_{\Lambda^a_\mu((F^{-1})^b)}\prec
\|f\|_{b^\alpha_{p,q}},
\end{equation}
where $a$, $b$ are defined by (\ref{eq412}) and $F$ is given by
(\ref{eq238}).
\end{corollary}

The results of Theorem \ref{thm52} and Corollary \ref{cor52}
for $q<p$ are not sharp in the sense that the range space can
be taken smaller in general. To see this, we are going to use a
slightly different approach.

We need  classes of measures $V^\Phi$, generated by the Lebesgue
measure in $\mathbf {R}^n,$ as follows:
$$
\mu\in V^\Phi \;\text{iff}\; \mu(TO)\prec
\Phi(|O|^{1/n})),
$$
uniformly for all open sets $O\subset \mathbf {R}^n$, where $|O|$ is
the Lebesgue measure of $O$.

Using the embedding
\begin{equation}\label{eq416}
b^\alpha_{p,q} \subset L^{r,q},\;   1/r=1/p-\alpha/n,\;
0<\alpha<n/p, \;1\leq p<\infty,\; 0<q\leq \infty
\end{equation}
where $L^{r,q}$ is the Lorentz space (see \cite{BL}), we see that
$$
|O|^{1-\alpha p/n}\prec cap(O;b^\alpha_{p,q}).
$$
In particular, if $\Phi$ satisfies (\ref{eq00}), then
\begin{equation}\label{eq417}
V^\Phi \subset Cap^{\Phi_\alpha}(b^\alpha_{p,q}) \subset C^\Phi,\;
\Phi_\alpha(s)=\Phi(s^{n-\alpha p}).
\end{equation}

\begin{theorem}\label{thm54}
Let $\phi$ satisfy (\ref{eq014}) and let $\Phi$ satisfy
(\ref{eq00}). If $\mu\in V^\Phi$, then
\begin{equation}\label{eq418}
\|f\ast\phi_t\|_{\Lambda^q_\mu((\Phi^{-1})^{n/r})}\prec
\|f\|_{L^{r,q}}, \;1<r<\infty,\; 0<q\leq \infty.
\end{equation}
Conversely, the estimate (\ref{eq418}) implies $\mu\in
V^\Phi.$
\end{theorem}

\begin{proof}
Starting with (\ref{eq47}), we get
$$
h(\lambda)\prec \Phi(|\{x: Mf(x)>c\lambda\}|^{1/n}).
$$
Hence
$$
\|f\ast\phi_t\|_{\Lambda^q_\mu((\Phi^{-1})^{n/r})}\prec
\|Mf\|_{L^{r,q}}\prec\|f\|_{L^{r,q}}.
$$

Choosing $f$ to be the characteristic function of the set $O$, we
derive as before that the estimate (\ref{eq418}) implies $\mu\in
V^\Phi.$ (see also \cite{S2}, Chapter II, Section 5.9, where the case
$\Phi(s)=s^{n\beta},$
$\beta>0$ is considered.)
\end{proof}

As a consequence we get
\begin{corollary}\label{cor53}
Let $\phi$ satisfy (\ref{eq014}) and let $\Phi$ satisfy
(\ref{eq00}) and let $\mu\in V^\Phi.$ Then
\begin{equation}\label{eq414}
\|f\ast\phi_t\|_{\Lambda^q_\mu((\Phi^{-1})^{q/p})}\prec
\|f\|_{b^\alpha_{p,q}}, \;1\leq p<n/\alpha,\; 0<q\leq \infty.
\end{equation}
\end{corollary}

\begin{corollary}\label{cor54}
Let $\phi$ satisfy (\ref{eq014}) and let $\Phi$ satisfy
(\ref{eq00}), (\ref{eq01}) and
\begin{equation}\label{eq015}
\sum \Phi(t_j^{1/n})\prec \Phi((\sum t_j)^{1/n}),\; t_j>0.
\end{equation}
If $\mu\in C^{\Phi_\alpha},$
$\Phi_\alpha(s)=\Phi(s^{n-\alpha p}),$ $1<p<n/\alpha,$ then
\begin{equation}\label{eq415}
\|f\ast\phi_t\|_{\Lambda^q_\mu((\Phi^{-1})^{q/p})}\prec
\|f\|_{b^\alpha_{p,q}}.
\end{equation}
\end{corollary}

Before proving this result, note that the function
$\Phi(s)=s^{n\beta} \log^\gamma (1+s),$ $\beta\geq 1,$ $\gamma\geq
0,$ satisfies the conditions (\ref{eq00}),  (\ref{eq01}),
(\ref{eq015}). Moreover, if $\beta>1,$ then we can take $\gamma$
to be any real number. Indeed, let us check (\ref{eq015}) for
$\beta>1.$ Choose $\epsilon>0$ such that $\beta-\epsilon/n>1,$ and
notice that the function $s^{\epsilon/n} \log^\gamma (1+s^{1/n})$
is equivalent to an increasing function, therefore
$$
t_j^{\epsilon/n} \log^\gamma(1+t_j^{1/n})\prec (\sum
t_j)^{\epsilon/n} \log^\gamma(1+(\sum t_j)^{1/n}),
$$
and then
$$
\sum t_j^{\beta-\epsilon/n} t_j^{\epsilon/n}
\log^\gamma(1+t_j^{1/n})\prec
\sum t_j^{\beta-\epsilon/n}(\sum t_j)^{\epsilon/n} \log^\gamma(1+(\sum
t_j)^{1/n})
$$
$$
\prec (\sum t_j)^{\beta} \log^\gamma(1+(\sum t_j)^{1/n}).
$$

To prove Corollary \ref{cor54}, we notice that as in \cite{S2}, we have
\begin{equation}\label{eq419}
V^\Phi=C^\Phi\;\text{if}\; \Phi \;\text{satisfies}\;
(\ref{eq00})\;\text{and}\; (\ref{eq015}).
\end{equation}
It remains to apply Theorem \ref{thm54}, the embeddings
(\ref{eq416}) and (\ref{eq417}), and the relation (\ref{eq419}).

\section{Negative results}

\begin{theorem}\label{thm61}
Let the continuous function $\Phi$ satisfy (\ref{eq00}) and let
\begin{equation}\label{eq016}
\int_0^1 \left[\frac{\Phi(u)}{u^{n-\alpha
p}}\right]^{\frac{1}{p-1}} \frac{du}{u}=\infty,\; 1<p<n/\alpha.
\end{equation}
Then
\begin{equation}\label{eq61}
Cap ^{\Phi_\alpha}(h^\alpha_p)\neq C^\Phi,\;
\Phi_\alpha(s)=\Phi(s^{\frac{1}{n-\alpha p}}).
\end{equation}
\end{theorem}

\begin{proof}
The condition (\ref{eq016}) and Theorem 5.4.2 of \cite{AH} imply the
existence of a compact $K\subset \mathbf {R}^n$ such that
$\Lambda_\Phi(K)>0$ and
\begin{equation}\label{eq62}
cap(K;h^\alpha_p)=cap(K;H^\alpha_p)=0.
\end{equation}
By Theorem 5.1.12 of \cite{AH}, there exists a positive measure
$\nu\in M^\Phi(h^\alpha_p)$, such that $\mu=\nu\times \delta_t \in
C^\Phi$ and $\nu(K)\approx \Lambda_\Phi(K).$ In particular,
\begin{equation}\label{eq63}
\nu(K)>0.
\end{equation}
Suppose that $\mu\in Cap ^{\Phi_\alpha}(h^\alpha_p).$ Then by
Theorem \ref{thm51}, where $\Phi$ is replaced by $\Phi_\alpha,$ we
have the estimate (\ref{eq46}). For any open set $O\supset K$, let
$f\in \cal S$ be such that $\|f\|^p_{h^\alpha_p} \leq 2
cap(O;h^\alpha_p).$ Then by (\ref{eq231}) and (\ref{eq46}),
$$
\nu(O)\prec h(d_\phi)\prec \Phi_\alpha(cap (O;h^\alpha_p)).
$$
Taking the monotone limit $O\mapsto K,$ we get $\nu(K)\leq 0,$
which contradicts (\ref{eq63}).
\end{proof}

For example, if $\Phi(u)=u^{n-\alpha p}(1+|\log u|)^\gamma,$
$\gamma+p-1\geq 0,$ $1<p<n/\alpha,$ then $\Phi_\alpha(u)\approx u
(1+|\log u|)^\gamma$ and $Cap^{\Phi_\alpha}(h^\alpha_p)\subset
C^\Phi,$ but these spaces are different.

\end{document}